\documentclass[11pt,a4paper]{article}
\usepackage{epsfig}
\usepackage[T1]{fontenc}    
\usepackage{graphics}
\usepackage{graphicx}
\usepackage{pstricks,pst-coil,pst-fill,pst-plot}
\usepackage[fleqn]{amsmath}    
\usepackage{amssymb}    
\usepackage{amsfonts}   
\usepackage{verbatim}   
\usepackage{mathrsfs}   
\usepackage{dsfont}
\usepackage{euscript}
\usepackage{yfonts}
\usepackage{enumerate}     
\usepackage{txfonts}
\usepackage{marvosym}
\usepackage{vmargin}        

\setmarginsrb{1.8cm}{2cm}{1.8cm}{2cm}{1cm}{1cm}{1cm}{1.6cm}
 \makeatletter
 \@addtoreset{equation}{section}
 \makeatother

\providecommand{\MR}{\relax\ifhmode\unskip\space\fi MR }

\providecommand{\href}[2]{#2}

\let\tend=\rightarrow

\long\def\symbolfootnote[#1]#2{\begingroup%
\def\thefootnote{\fnsymbol{footnote}}\footnote[#1]{#2}\endgroup}

\newtheorem{theorem}{Theorem}[section]
\newtheorem{prop}{Proposition}[section]
\newtheorem{cor}{Corollary}[section]
\newtheorem{defin}{Definition}[section]

\newtheorem{lemme}{Lemma}[section]

\def\Proof{\medskip\noindent {\it Proof --- \ }}

\def\qed{\hfill\rule{2mm}{2mm}}

\newcommand\beq{\begin{equation}}
\newcommand\enq{\end{equation}}
\newcommand\bem{\begin{multline}}
\newcommand\enm{\end{multline}}

\def\beqa{\begin{eqnarray}}
\def\eeqa{\end{eqnarray}}
\def\ba{\begin{array}}
\def\ea{\end{array}}
\def\det{\operatorname{det}}

\newcommand{\f}[2]{{\ensuremath{%
    \mathchoice%
    {\dfrac{#1}{#2}}
    {\dfrac{#1}{#2}}
    {\frac{#1}{#2}}
    {\frac{#1}{#2}}
}}}
\newcommand{\tf}[2]{\ensuremath{#1/#2}}

\def\a{\alpha}

\def\Ga{\Gamma}

\def\de{\delta}

\def\eps{\epsilon}

\def\la{\lambda}

\def\sg{\sigma}

\def\Ups{\Upsilon}

\newcommand{\mc}[1]{\ensuremath{\mathcal{#1}}}

\newcommand{\msc}[1]{\ensuremath{\mathscr{#1}}}

\newcommand{\bs}[1]{\ensuremath{\boldsymbol{#1}}}

\newcommand{\ov}[1]{\ensuremath{\overline{#1}}}
\newcommand{\wt}[1]{\ensuremath{\widetilde{#1}}}
\newcommand{\wh}[1]{\ensuremath{\widehat{#1}}}

\newcommand{\Int}[2]{\ensuremath{\int\limits_{#1}^{#2}}}
\newcommand{\Oint}[2]{\ensuremath{\oint\limits_{#1}^{#2}}}

\newcommand{\sul}[2]{\ensuremath{\sum\limits_{#1}^{#2}}}

\newcommand{\Cx}{\ensuremath{\mathbb{C}}}

\newcommand{\Dp}[1]{\ensuremath{\partial_{#1}}}

\newcommand{\ex}[1]{\ensuremath{\e{e}^{#1}}}

\newcommand{\norm}[1]{\ensuremath{|| #1 ||}}

\newcommand{\dd}{\mathrm{d}}
\newcommand{\e}[1]{\ensuremath{\mathrm{#1}}}

\newcommand{\intn}[2]{\ensuremath{[\![ \, #1 \,;\, #2 \,]\!]}}

\begin{document}

\begin{flushright}

\end{flushright}
\par \vskip .1in \noindent

\vspace{14pt}

\begin{center}
\begin{LARGE}
{\bf On lacunary Toeplitz determinants.}
\end{LARGE}

\vspace{30pt}

\begin{large}

{\bf K.~K.~Kozlowski}\footnote[1]{Universit\'{e} de Bourgogne, Institut de Math\'{e}matiques de Bourgogne, UMR 5584 du CNRS, France,
karol.kozlowski@u-bourgogne.fr}. 
\par

\end{large}

\vspace{40pt}

\centerline{\bf Abstract} \vspace{1cm}
\parbox{12cm}{\small
By using Riemann--Hilbert problem based techniques, we obtain the asymptotic expansion of lacunary Toeplitz determinants 
$\det_N\big[ c_{\ell_a-m_b }[f] \big]$  
generated by holomorhpic symbols, where $\ell_a=a$ (resp. $m_b=b$) except for a finite subset of indices $a=h_1,\dots, h_n$ 
(resp. $b=t_1,\dots, t_r$). 
In addition to the usual Szeg\"{o} asymptotics, our answer involves a determinant of size $n+r$.}

\end{center}

\vspace{40pt}

\section*{Introduction}

A lacunary Toeplitz determinant generated by a symbol $f$ refers to the below determinant
\beq
\det_{N}\Big[ c_{\ell_a-m_b}[f] \Big] \qquad \e{where} \qquad 
c_n[f] \; = \; \Oint{ \Dp{}\mc{D}_1 }{} \f{ f(z)}{ z^{n+1} } \cdot  \f{ \dd z  }{2i\pi} 
\label{ecriture detereminant lacunaire general}
\enq
and $\Dp{}\mc{D}_{\eta}$ is the counter clockwise oriented boundary of the disc of radius $\eta$
centred at $0$. 
The sequences $\ell_a$, $m_b$ appearing in \eqref{ecriture detereminant lacunaire general} are such that 
\beqa
\ell_a & = & a \quad \e{for} \quad a\in \big\{ 1, \dots, N \big\} \setminus \big\{h_1,\dots, h_n \big\} 
\qquad \e{and} \qquad 
\ell_{h_a} \; = \; p_a \quad a=1,\dots, n \;  \label{definition de la suite ella}\\
m_a & = & a \quad \e{for} \quad a\in \big\{ 1, \dots, N \big\} \setminus \big\{t_1,\dots, t_r \big\} 
\qquad \e{and} \qquad 
m_{t_a} \; = \; k_a \quad a=1,\dots, r \; \label{definition de la suite ma}
\eeqa
The integers $h_a \in \intn{1}{N}$ and $p_a \in \mathbb{Z} \setminus \intn{1}{N}$, $a=1,\dots, n$ 
(resp. $t_a \in \intn{1}{N}$ and $k_a \in \mathbb{Z} \setminus \intn{1}{N}$, $a=1,\dots, r$)
are assumed to be pairwise distinct. 
The large-$N$ asymptotic behaviour of such determinants has been first considered 
by Tracy and Widom \cite{TracyWidomAsymptoticExpansionLacunaryToeplitz} and Bump and Diaconis \cite{BumpDiaconisLacunaryToeplitzThrougSumsSymFctsAndYoungTableaux}. More or less at the same time,
these authors have obtained two formulae of a very different kind for these large-N asymptotics. 
In fact, both collaborations expressed the large-$N$ behaviour of the lacunary Toeplitz determinant
in terms of the unperturbed determinant $\det_{N}\Big[ c_{a - b}[f] \Big]$  times 
an extra term whose representations took a very different form. 
The expression found  Bump and Diaconis was based on characters of the symmetric group 
associated with the partitions $ \la $ and $\mu$ that can be naturally associated 
with the sequences $\ell_a$ and $m_b$. The answer involved the sum over the symmetric groups of $|\la|$
and $|\mu|$ elements. In their turn, Tracy and Widom obtained a determinant representation of the type
\beq
\det_{N}\big[ c_{\ell_a-b}[f] \big] \; = \; \det_{N-m}\big[ c_{j-k}[f] \big] \cdot 
\det_{ q }\big[ W_{jk} \big] \cdot  \Big(  1+ \e{o}(1) \Big)    \qquad q\, = \, \max \big\{ t_1,\dots , t_r , h_1,\dots, h_n   \big\} 
\label{formule asymptotique Widom}
\enq
 where $W_{jk} $ was an explicit $q \times q$ sized matrix depending on the symbol $f$ and the numbers
$h_1, \dots, h_n$, $p_1,\dots, p_n$, $t_1,\dots, t_r$ and $k_1,\dots, k_r$. In \cite{DehayeProofIdentityBumpDiaconisTracyWidomLacunatyToeplitz},  Dehaye
proved, by a direct method, the equivalence between the two aforementioned formulae. One should also mention that the 
large-$N$ asymptotic behaviour of some generalizations of lacunary Toeplitz determinants have been obtained by Lions 
in \cite{LionsToeplitzLacunaires}. 

The drawbacks of the aforementioned asymptotic expansions was that the answer 
depended on the magnitude of the lacunary parameters $p_a, k_b, h_a, t_b$. As soon as these
parameters were also growing with $N$, the form of the answer did not allow for an easy access to the large-$N$
asymptotic behaviour of the lacunary determinant. Indeed, in Bump-Diaconis' case, the number of summed up terms was
growing as $\sum(p_a-h_a) + \sum(k_a-t_a)$ whereas in Tracy-Widom's case, the non-trivial determinant part 
involved a matrix of size $\max\{h_a, t_b\}$. 

In the present note we obtain an asymptotic expansion solely in terms of a $(n+r)\times (n+r)$ matrix
and show that the latter is enough so as to treat certain cases of lacunary parameters $p_a, k_b, h_a, t_b$ going to infinity. 
The structure of the asymptotics when $r\not=0$ (\textit{ie} $m_a\not=a$) is slightly more complex, so that 
we postpone the statement of the corresponding results to the core of the paper and present the asymptotic
expansion we obtain on the example of line-lacunary Toeplitz determinants 
\begin{theorem}
\label{Theorem Asymptotiques Toeplitz lacunaire a lignes}
Let $f$ be a non-vanishing function on $\Dp{}\mc{D}_{1}$ such that $f$ and $\ln f$
are holomorphic on some open neighbourhood of $\Dp{}\mc{D}_{1}$. 
Let $\ell_a$ be defined as \eqref{definition de la suite ella} and $\a$
be the piecewise analytic function 
\beq
\a(z) \; = \;   \exp\bigg\{ -  \sul{n \geq 0}{} c_{ n}\big[ \ln f \big] \cdot z^{n}  \bigg\} 
\quad  \e{for}\;  z \in \mc{D}_{1}  \qquad \e{and} \qquad 
			\a(z) \; = \;	  \exp\bigg\{  \sul{n \geq 1}{} c_{- n}\big[ \ln f \big] \cdot z^{-n}  \bigg\}  
								  \quad  \e{for}\; z \in \Cx\setminus \ov{\mc{D}}_{1}   \;  . 
\label{definition facteur alpha solution RHP scalaire}
\enq
Then, provided that the matrix $M$ given below is non-singular, the line-lacunary Toepltz determinant 
$\det_{N}\big[ c_{\ell_a-b}[f] \big]$ admits the representation 
\beq
\det_{N}\big[ c_{\ell_a-b}[f] \big] \; = \; \det_{N}\big[ c_{a-b}[f] \big] \cdot  \det_{n}\big[ M_{ab} \big] \cdot  
\Big(  1+ \e{O}\big( N^{-\infty} \big) \Big) \;, 
\label{Theorem intro ecriture forme DA}
\enq
where the $n\times n$ matrix $M$ reads 
\bem
M_{ab} \; = \; - \bs{1}_{\mathbb{N}}(p_a) 
\Oint{  \Dp{}\mc{D}_{\eta_z}    }{} \hspace{-2mm} \f{ \dd z }{ 2i\pi} \cdot 
\Oint{  \Dp{}\mc{D}_{\eta_s}    }{} \hspace{-2mm} \f{ \dd s }{ 2i\pi} 
\cdot \f{ \a(z) }{ \a(s) } \cdot \f{ s^{N-p_a} \cdot z^{h_b-N-1}  }{ z- s }    \\
\; + \; \bs{1}_{\mathbb{N}}(-p_a) \Oint{  \Dp{}\mc{D}_{ \eta_z^{-1} }    }{} \hspace{-2mm} \f{ \dd z }{ 2i\pi} \cdot 
\Oint{  \Dp{}\mc{D}_{ \eta_s^{-1} }    }{} \hspace{-2mm} \f{ \dd s }{ 2i\pi} 
\cdot \f{ \a(s) }{ \a(z) } \cdot \f{ s^{-p_a} \cdot z^{h_b-1}  }{ z- s }  \;, 
\label{ecriture formule asymptotique matrice M cas lacunaire a ligne}
\end{multline}
 and $1>\eta_{z} > \eta_{s} >0 $\;. 

\end{theorem}

The theorem above allows one to obtain the large $N$-asymptotic expansion of the line-lacunary Toeplitz determinant
independently on the magnitude (in respect to $N$) of the lacunary parameters $\{h_a\}$ and $\{ p_a \}$. 
Indeed, since the size of the matrix $M$ does not depend on the integers $\{h_a\}$ or $\{p_a\}$, the problem boils down to a \textit{classical}
asymptotic analysis of one-dimensional integrals. Still, in order to provide one with an explicit answer, 
some more data on these parameters is needed. For instance, one has the 
\begin{cor}
\label{Corolaire Toeplitz lacunaire à un indice}
Let
\beqa
p_a   \; = \; 1-p_a^{-} \quad a=1,\dots, n_{-} \qquad &\e{and}& \qquad p_{a+n_-} \; = \;  p_a^{+} + N \quad a=1,\dots, n_{+}  
\label{particules trous pour entiers pa} \\
h_a   \; = \; h_a^{-} \quad a=1,\dots, n_{-} \qquad &\e{and} &\qquad h_{a+n_-} \; = \; N+1- h_a^{+} \quad a=1,\dots, n_{+}   \;, 
\label{particuler trous pour entiers ha}
\eeqa
where $p_a^{\pm}$ and $h_a^{\pm}$ are assumed to be independent of $N$ and $n=n_- + n_+$. 
Provided that the matrices $M^{(\pm)}$ given below are not singular, one has 

\beq
 \det_{n}\big[ M_{ab} \big] \; = \;\det_{n_+}\Big[ M_{ a  b }^{(+)} \Big]  
  \cdot \det_{n_-}\Big[ M_{ a  b }^{(-)} \Big] \cdot \Big( 1+ \e{O}\big( N^{-\infty} \big) \Big) \;, 
\enq 
where 
\beq
 M_{a b }^{(+)} \; = \; - 
 \Oint{   \Dp{}\mc{D}_{\eta_z}  }{} \hspace{-2mm}\f{ \dd z }{2i\pi}  \hspace{-2mm} \cdot 
 \Oint{   \Dp{}\mc{D}_{\eta_s }  }{} \hspace{-2mm} \f{ \dd s }{2i\pi} \cdot 
\f{ s^{ -p_{a}^+ }\cdot z^{- h_{b}^+} }{  z \, - \,  s }     \cdot 
\f{ \a(z) }{ \a(s) }   \qquad \e{and} \qquad 
M_{a b }^{(-)} \; = \; 
 \Oint{   \Dp{}\mc{D}_{\eta_z^{-1}}  }{} \hspace{-2mm}\f{ \dd z }{2i\pi}  \hspace{-2mm} \cdot 
 \Oint{   \Dp{}\mc{D}_{\eta_s^{-1} }  }{} \hspace{-2mm} \f{ \dd s }{2i\pi} \cdot 
\f{ s^{p_{a}^- - 1}\cdot z^{ h_{b}^- - 1} }{  z \, - \,  s }     \cdot 
\f{ \a(s) }{ \a(z) }  \;. 
\enq
\end{cor}

We obtain the asymptotic expansion \eqref{Theorem intro ecriture forme DA} 
by interpreting the lacunary Topelitz determinant as the determinant of 
a finite rank perturbation of a integrable integral operator acting on the unit circle. 
The inverse of the integrable integral operator, in the large-$N$ regime, can be constructed by means of 
an asymptotic resolution of a Riemann--Hilbert problem. We have restricted the study of the present paper 
to holomorphic symbols. However, in principle, one could apply the method to less regular symbols, \textit{eg}
those containing Fischer-Hartwig singularities. Of course, the price of such generalisation would be to deal 
with certain technicalities related with the more complex structure of the large-$N$ approximant to the associated
resolvent operator. 

The paper is organized as follows. We prove Theorem \ref{Theorem Asymptotiques Toeplitz lacunaire a lignes}  
in Section \ref{Section Toeplitz lacunaire a lignes}. In Section \ref{Section Toeplitz lacunaire a lignes et colonnes}
we establish the large-$N$ asymptotic expansion of general line and row lacunary Toeplitz determinants 
subordinate to the sequences \eqref{definition de la suite ella}-\eqref{definition de la suite ma}. 
Technical details related to the large-$N$ inversion of integrable integral operators 
arising in the analysis of Toeplitz determinant generated by holomorphic non-vanishing on $\Dp{}\mc{D}_1$ symbols are 
recalled in appendix \ref{Appendix RHP pour Toeplitz regulier}.

\section{The line lacunary Toeplitz determinants}
\label{Section Toeplitz lacunaire a lignes}

In this section, we first prove a preliminary factorisation result that allows one to express the lacunary Toeplitz determinant 
$\det_{N}\Big[  c_{ \ell_a - b }[f] \Big]$ in terms of the non-perturbed Toeplitz determinant  $\det_{N}\Big[  c_{ a - b }[f] \Big] $ 
and of the determinant of a $n\times n$ matrix. We subsequently analyse the large-$N$ behaviour of this finite-size $n$
determinant. 

\subsection{The factorisation}

\begin{lemme}
\label{Lemme factorisation Toeplitz lacunaire a lignes}
Let $f$ be non-vanishing on $\Dp{}\mc{D}_1$ and such that $f$ and $\ln f$ are holomorphic in some open neighbourhood of $\msc{C}$.
Let $V_0$ be the integral kernel 
\beq
V_0\big( z, s \big)  \; = \; \big( f(z) - 1 \big) \cdot 
\f{  z^{\f{N}{2}}\cdot s^{-\f{N}{2}} \; - \; z^{-\f{N}{2}}\cdot s^{\f{N}{2}}  }
{ 2i\pi \big( z - s \big) }
\label{definition noyau integral V0}
\enq
of the integrable integral operator $I+V_0$ acting on $L^2\big( \Dp{}\mc{D}_1 \big)$. 
Then, provided that $N$ is large enough, $I+V_0$ is invertible with inverse $I-R_0$ and the below
factorization holds
\beq
\det_{N}\Big[  c_{ \ell_a - b }[f] \Big] \; = \; \det_{N}\Big[  c_{ a - b }[f] \Big] \cdot 
\det_n \big[ M_{ab} \big]
\enq
where 
\beq
M_{ k \ell } \; = \; \de_{k \ell } \; - \; c_{ h_{k} - h_{\ell} }[f] \; + \; c_{ p_{k} - h_{\ell} }[f]
\; + \;  \Int{ \msc{C} }{} R_{0}(z,s) \cdot f(s)\cdot \big( s^{\f{N}{2}-h_k} - s^{\f{N}{2} - p_k} \big) \cdot z^{h_{\ell}- 1 -\f{N}{2}} \cdot 
\f{ \dd s \cdot \dd z  }{ 2i \pi } \;. 
\label{definition matrice Mkl}
\enq
\end{lemme}

\Proof

Let $I+V$ be the integral operator on $L^{2}(\Dp{}\mc{D}_{1})$ with a kernel given by 
\beq
V(z,s) \; = \; \sul{a=1}{N} \kappa_a(z) \cdot \tau_a(s)  \qquad \e{where} \quad  \quad 
\tau_a(z) \; = \; \f{1}{2i\pi} \cdot z^{a-1-\f{N}{2}}
\enq
and
\beq
\ba{cccc}  \kappa_a(z) & = & \big( f(z)-1 \big) \cdot z^{\f{N}{2}-a}   
											& a \in \{1,\dots , N \}  \setminus \{h_1, \dots, h_n \} \hspace{3mm} \\
			\kappa_{h_a}(z) & = & f(z) \cdot  z^{ \f{N}{2}-p_a } \; - \; z^{ \f{N}{2}-h_a }    & a=1, \dots, n  	 \ea    \;. 
\enq
Since $V$ is a finite rank $N$ operator, the Fredholm determinant of $I+V$ reduces to one of an $N\times N$ matrix
\beq
\det_{\Dp{}\mc{D}_{1}}\big[  I \, + \, V \big] \; = \; 
\det_{N}\Big[ \de_{ab} \; + \; \int_{ \Dp{}\mc{D}_{1} }{} \kappa_{a}(z) \cdot  \tau_b(z) \cdot \dd z \Big]
\; = \; \det_{N}\Big[  c_{ \ell_a - b }[f] \Big] \;. 
\enq
One can decompose the kernel $V$ as $V=V_0 + V_1$ where $V_0$ has been introduced in \eqref{definition noyau integral V0}
whereas $V_1$ is the finite rank $n$ perturbation of $V_0$ given by 
\beq
V_1\big( z, s \big)  \; = \;  -\f{ f(z) }{ 2i\pi }  \cdot \sul{ a =1 }{ n } 
\big( z^{\f{N}{2} -h_a} \; - \;  z^{\f{N}{2} -p_a}  \big) 
\cdot s^{h_a -1 -\f{N}{2} }  \;. 
\enq
It follows from the strong Szeg\"{o} limit theorem and from the identity 
$\det_{ \Dp{}\mc{D}_{1} }\big[ I\, + \, V_0 \big]  \; = \; \det_{N}\Big[  c_{ a - b }[f] \Big]$
that, provided $N$ is taken large enough, the operator $I+V_0$ is invertible. Hence, all-in-all, we get that
\beq
\det_{ \Dp{}\mc{D}_{1} }\big[ I \, + \, V \big] \; = \;   \det_{ \Dp{}\mc{D}_{1} }\big[ I\, + \, V_0 \big] 
\cdot \det_{ \Dp{}\mc{D}_{1} }\big[ I\, + \, (I-R_0)\cdot V_1 \big] 
\; = \;  \det_{N}\Big[  c_{ a - b }[f] \Big]  \cdot \det_n\big[ M_{ k \ell } \big] 
\enq
where the matrix $M_{k\ell}$ is as defined in \eqref{definition matrice Mkl}.

\subsection{Asymptotic analysis of $\det_{n}[M]$-Proof of theorem \ref{Theorem Asymptotiques Toeplitz lacunaire a lignes}}

As it has been recalled in the appendix, the resolvent kernel $R_0$ of the operator $I+V_0$ can be recast as 
\beq
R_0 \; = \; R_0^{(0)} + R_0^{(\infty)}
\label{ecriture decomposition resolvant en partie finie et perturbation exp petite en N}
\enq
where 
\beq
R_0^{(0)}(z,s) \; = \;  \f{ f(z) -1 }{ 2i\pi } \cdot 
\f{ z^{\f{N}{2}} \cdot s^{-\f{N}{2}} \cdot  \a_+(s) \cdot  \a_-^{-1}(z) 
							\; - \;   
							s^{\f{N}{2}} \cdot  z^{-\f{N}{2}} \cdot  \a_+(z) \cdot \a_-^{-1}(s) }{ z-s }
\label{definition resolvent approximatif de V0}
\enq
and
\beq
\norm{ R_0^{(\infty)} }_{L^{\infty}\big( \Dp{}\mc{D}_{1} \times \Dp{}\mc{D}_{1} \big) } \; \leq \; C \cdot N \cdot  \ex{-\kappa N } \;. 
\enq
Above, $\a_{\pm}$ are the $+$ (\textit{ie} from within) and $-$ (\textit{ie} from the outside) non-tangential 
limits on $\Dp{}\mc{D}_1$ of the piecewise analytic function $\a$ defined in \eqref{definition facteur alpha solution RHP scalaire}.  
The decomposition \eqref{ecriture decomposition resolvant en partie finie et perturbation exp petite en N} ensures that, 
for the price of exponentially small corrections, one can trade the kernel $R_0$ for $R_0^{(0)}$ in \eqref{definition matrice Mkl}. 
Using that $\a_{+}$  (resp. $\a_-$) admit an analytic continuation to some open neighbourhood of $\Dp{} \mc{D}_1$ 
in $\Cx \setminus \ov{\mc{D}}_1$ (resp. interior $\mc{D}_1$) we deform the contours in the double integral
associated with $R_0^{(0)}$ to 
\begin{itemize}
\item $ \Dp{}\mc{D}_{\eta_z^{-1}} \times  \Dp{}\mc{D}_{\eta_s^{-1}}$ in what concerns the part of the integrand 
containing $\a_+(s)/\a_-(z)$; 
\item  $ \Dp{}\mc{D}_{\eta_z} \times  \Dp{}\mc{D}_{\eta_s}$in what concerns the part of the integrand 
containing $\a_+(z)/\a_-(s)$.  
\end{itemize}
 The resulting residue cancels out the pre-factors in \eqref{definition matrice Mkl}
leading to 
\bem
M_{k\ell} \; = \;  \Oint{  \Dp{}\mc{D}_{\eta_s^{-1}} }{}\hspace{-2mm} \f{\dd s }{2i \pi } 
 \Oint{  \Dp{}\mc{D}_{\eta_z^{-1}} }{} \hspace{-2mm} \f{\dd z }{2i \pi }
\; \a_-(s) \cdot \big( \a_+^{-1}(z) \, - \, \a_-^{-1}(z)\big)\cdot \big(s^{-h_k} \, -  \,  s^{-p_k} \big) \cdot 
\f{ z^{h_{\ell}-1} }{z-s}  \\
\; +  \; \Oint{  \Dp{}\mc{D}_{\eta_s} }{}\hspace{-2mm} \f{\dd s }{2i \pi } 
 \Oint{  \Dp{}\mc{D}_{\eta_z} }{} \hspace{-2mm} \f{\dd z }{2i \pi }
\; \a_+(s)^{-1} \cdot \big( \a_+(z) \, - \, \a_-(z)\big)\cdot \big(s^{N-h_k} \, -  \,  s^{N-p_k} \big) \cdot 
\f{ z^{h_{\ell}-1-N} }{z-s} \; \; + \; \; \e{O}(N^{-\infty}) \;. 
\label{ecriture formule intermediaire pour Mab}
\end{multline}
The term $s^{-h_k}$ (resp. $s^{N-h_k}$) do not contribute to the integral as can be seen by deforming the contour of 
$s$-integration to $\eta_s$ (resp. $\eta_s^{-1}=0$). Further, the first line of \eqref{ecriture formule intermediaire pour Mab}
only gives non-vanishing contributions if $p_k \leq 0$ (resp. the last line of \eqref{ecriture formule intermediaire pour Mab}
only gives non-vanishing contributions if $p_k \geq N+1$). 
This yields \eqref{ecriture formule asymptotique matrice M cas lacunaire a ligne}. \qed






\section{The asymptotic expansion of line and row lacunary Toeplitz determinants}
\label{Section Toeplitz lacunaire a lignes et colonnes}

\subsection{The factorisation in the general case}

The factorized representation in the general case depends, in particular, on whether there are some overlaps 
between the integers parametrising the lacunary line and columns. 
We thus need a definition so as to be able to distinguish between the different cases. 

\begin{defin}
The sets $ \{ h_a \}_1^{n}$ and $\{ t_b \}_1^r$  with $h_a, t_b \in \intn{1}{N}$
are said to be well-ordered with overlap $c \in \intn{0}{\min\{r,n\} }$ if  
\beq
 h_a \; = \; t_a \quad \e{for} \quad  a=1,\dots, c  \qquad \e{whereas}  \qquad
 \big\{ h_{c+1},\dots, h_n \big\} \cap  \big\{ t_{c+1},\dots, t_r \big\} \; = \; \emptyset \;. 
\label{definition parametre overlap c}
\enq

\end{defin}

It is clear that given two not well ordered sequences $\ell_a$ and $m_a$, one can always relabel
the indices of the lacunary integers $\{p_a, h_a, k_b, t_b\}$ so that \eqref{definition parametre overlap c}
holds. 
There is thus no restriction in assuming that the sequences $\ell_a$ and $m_a$ are well ordered, so that 
we are going to do so in the following.

\begin{prop}
\label{Proposition factorisation Toeplitz lacunaire a ligne et colonnes}

Let $\ell_a$ and $m_a$ be sequences as defined in \eqref{definition de la suite ella}-\eqref{definition de la suite ma}
and $f$ a non-vanishing holomorphic function on some open neighbourhood of $\Dp{}\mc{D}_{1}$ 
such that $\ln f $ is also holomorphic on this neighborhood. Then, the lacunary line and row Toeplitz determinant 
admits the representation 
\beq
\det_{N} \Big[ c_{\ell_a - m_b }[f]  \Big] \; = \; \det_{N}\Big[  c_{a - b }[f]  \Big] \cdot 
\det_{n+r}\big[  \mc{N}  \big] \;. 
\label{ecriture factorisation explicite Toeplitz lacunaire ligne et colonnes}
\enq
The matrix $\mc{N}$ appearing above admits the blocks structure
\beq
\mc{N} \; = \; \left(   \ba{c c}   \mc{N}_{I ; I}  & \mc{N}_{I ; II}  \\ 
									\mc{N}_{II ; I}  & \mc{N}_{II ; II}	\ea	\right) 
\enq
with blocks being given by 
\beq
\big( \mc{N}_{A ; I} \big)_{ab} \; = \; \de_{A;I} \de_{ab} \de_{b>c} \; + \; \Oint{ \Dp{}\mc{D}_1 }{} 
U_{A;a}(z) \cdot v_{I;b}(z) \cdot \dd z \quad \e{and} \quad 
\big( \mc{N}_{A ; II} \big)_{ab} \; = \; \de_{A;II} \de_{ab} \de_{b \leq c} \; + \; \Oint{ \Dp{}\mc{D}_1 }{} 
U_{A;a}(z) \cdot v_{II;b}(z) \cdot \dd z  
\label{ecriture entrees matrice N}
\enq
in which $A \in \{ I , II \}$. The functions $U_{A;a}$ are built in terms of the resolvent $R_0$ 
to the integral operator $I+V_0$ defined in \eqref{definition noyau integral V0}
and of the functions $u_{A;a}$
\beq
u_{I;a}(z) \; = \; \f{ f(z) }{ 2i\pi } \cdot z^{\f{N}{2} -p_a} \; - \; 
 \f{ z^{\f{N}{2}-h_{a}} }{ 2i\pi } \cdot \Big(\de_{a\leq c} \;+ \; \de_{a>c} f(z) \Big) \qquad 
\e{and} \qquad 
u_{II;a}(z)\; = \;  \f{ f(z)-1 }{ 2i\pi } \cdot  z^{\f{N}{2} - t_a} 
\label{definition des fonctions u}
\enq
as $U_{A;a} \; = \;  (I-R_0)[u_{A;a}]$. Finally, the function $v_{A;b}$ appearing in \eqref{ecriture entrees matrice N} read
\beq
v_{I;b}(z) \; = \; \de_{b \leq c } \cdot z^{k_b - \f{N}{2} -1} \; + \; \de_{b>c}  \cdot z^{h_{b}-1-\f{N}{2}} \qquad 
\e{and} \qquad 
v_{II;b}(z)\; = \; - \de_{b \leq c } \cdot  z^{t_b - \f{N}{2} -1} \; + \; \de_{b>c} \cdot z^{k_{b}-1-\f{N}{2}} \; . 
\label{definition des fonction v}
\enq

\end{prop}

\Proof 

Let $I+\wh{V}$ be the integral operator acting on $L^{2}(\Dp{}\mc{D}_{1})$ with the kernel 
\beq
\wh{V}(z,s) \; = \; \sul{a=1}{N} \wh{\kappa}_a(z)\cdot \wh{\tau}_a(s)  \qquad \e{where} \quad  \left\{ \ba{cc} 
			\wh{\tau}_a(z) \; = \; \tf{ z^{a-1-\f{N}{2}} }{(2i\pi)}  & a \in \{1,\dots , N \}  \setminus \{t_1, \dots, t_r \} \vspace{2mm} \\
			\wh{\tau}_{t_a} (z) \; = \; \tf{ z^{k_a-1-\f{N}{2}} }{(2i\pi)}  & a = 1, \dots, r \ea \right. 
\enq
and
\beq
\ba{cccc}  \wh{\kappa}_a(z) & = & \big( f(z)-1 \big) \cdot z^{\f{N}{2}-a}  \; - \; {\displaystyle \sul{s=1}{r}} \de_{a,t_{s}} z^{\f{N}{2}-k_{s}}
					 & a \in \{1,\dots , N \}  \setminus \{h_1, \dots, h_n \} \vspace{1mm} \\
			\wh{\kappa}_{h_a}(z) & = &   f(z) \cdot  z^{ \f{N}{2}-p_a } \; -  \; z^{ \f{N}{2}-h_a } 
										\; - \; {\displaystyle \sul{s=1}{r}} \de_{h_a,t_{s}} z^{\f{N}{2}-k_{s}}	 & a=1, \dots, n  	 \ea    \;. 
\enq
It is readily seen that 
\beq
\det_{ \Dp{}\mc{D}_{1} }\big[  I \, + \, \wh{V} \big] \; = \; 
\det_{N}\Big[ \de_{ a b } \; + \; \int_{ \Dp{}\mc{D}_{1} }{} \wh{\kappa}_{a}(z) \cdot \wh{\tau}_{b}(z) \cdot \dd z \Big]
\; = \; \det_{N}\Big[  c_{ \ell_a  - m_b}[f] \Big] \;. 
\enq
The kernel $\wh{V}$ can be recast as $\wh{V}=V_0 + \wh{V}_1$ where $V_0$ has been introduced in \eqref{definition noyau integral V0}
and $\wh{V}_1$ is the finite rank $n+2r$ perturbation of $V_0$ given by 
\beq
\wh{V}_1\big( z, s \big)  \; = \;  \sul{ a =1 }{ n } u_{I;a}(z) \cdot v_{I;a}(s)
\; + \; \sul{ a =1 }{ r } \Big\{ u_{II;a}(z) \cdot \wt{v}_{II;a}(s) \; + \;  u_{III;a}(z) \cdot \wt{v}_{III;a}(s) \Big\} \;. 
\enq
The functions $u_{I,a}$, $u_{II,a}$ and $v_{I,a}$ are as defined in 
\eqref{definition des fonctions u}-\eqref{definition des fonction v} whereas
\beq
\wt{v}_{II;b}(z) \; = \; \de_{b>c}\cdot z^{k_b-\f{N}{2}-1} \; - \;   z^{t_b-\f{N}{2}-1}
\quad , \quad 
\wt{v}_{III;b}(z) \; = \; z^{k_b-\f{N}{2}-1}
\quad \e{and} \quad 
\wt{u}_{III;a}(z) \; = \;  - \f{ z^{\f{N}{2}-k_a} }{ 2i \pi }  \;. 
\enq
Proceeding as in the proof of lemma \ref{Lemme factorisation Toeplitz lacunaire a lignes}, we get that
\beq
\det_{ \Dp{}\mc{D}_{1} }\big[ I \, + \, V \big] \; = \;   \det_{ \Dp{}\mc{D}_{1} }\big[ I\, + \, V_0 \big] 
\cdot \det_{ \Dp{}\mc{D}_{1} }\big[ I\, + \, (I-R_0)\cdot \wh{V}_1 \big] 
\; = \;  \det_{N}\Big[  c_{ j-k }[b] \Big]  \cdot \det_{n+2r}\big[ \mc{M} \big] 
\enq
where $\mc{M}$ is the $(n+2r)\times (n+2r)$ block matrix
\beq
\mc{M} \; = \; \left(   \ba{c c c}   \mc{M}_{I ; I}  & \mc{M}_{I ; II} & \mc{M}_{I; III} \\ 
									\mc{M}_{II ; I}  & \mc{M}_{II ; II}	& \mc{M}_{II ; III}  \\ 
									\mc{M}_{III ; I}  & \mc{M}_{III ; II}	& \mc{M}_{III ; III}   \ea	\right) 
\quad \e{with} \quad 
\big( \mc{M}_{A;B} \big)_{ab} = \de_{A,B} \de_{ab} \; + \; \Oint{ \Dp{} \mc{D}_1 }{} U_{A;a}(z) \cdot v_{B;b}(z) \cdot \dd z	\;	.
\enq
The upper-case entries $A,B$ belong to $\{I, II, III\}$ whereas the lower-case entries $a,b$ subordinate to 
the upper-case entry $I$ run from $1$ to $n$ and those subordinate to the upper-case entries $II$ or $III$ run from $1$ to $r$. 
Since $\wt{u}_{III,a} \in \ker(V_0)$, it follows that $\wt{u}_{III,a} \in \ker(R_0)$. Then, a straightforward calculation shows that,
in fact, the block structure of the lines of type III simplifies leading to 
\beq
\mc{M} \; = \; \left(   \ba{c c c}   \mc{M}_{I ; I}  & \mc{M}_{I ; II} & \mc{M}_{I; III} \\ 
									\mc{M}_{II ; I}  & \mc{M}_{II ; II}	& \mc{M}_{II ; III}  \\ 
									\ba{cc} -I_{c} & 0 \\ 0 & 0  \ea  & \ba{cc} 0 & 0 \\ 0 & -I_{r-c}  \ea 	& 0   \ea	\right) \;. 
\enq
There, $I_k$ refers to the identity matrix in $k$ dimensions. In order to reduce the size of the determinant of $\mc{M}$
it is enough to exchange the first $c$ columns of the block $I$  with the first $c$ ones of the block $III$, and 
then exchange the $r-c$ last columns of block $II$ with the $r-c$ last columns of the block $III$. This produces, all in all, 
an overall $(-1)^r$ sign. The latter cancels out with the one issuing from $\det_r[-I_r]$, thus 
leading to $\det_{n+2r}[\mc{M}] = \det_{n+r}[\mc{N}]$. \qed

\subsection{A special case of a large $N$ asymptotics}

Replacing the exact resolvent $R_0$ by its approximate resolvent \eqref{definition resolvent approximatif de V0} 
in the definition of the matrix entries of $\mc{N}$ \eqref{ecriture factorisation explicite Toeplitz lacunaire ligne et colonnes} 
leads to exponentially small in $N$ corrections. 
Upon such a replacement, Proposition \ref{Proposition factorisation Toeplitz lacunaire a ligne et colonnes}
basically yields the most general expression for the large-$N$ asymptotics of the ratio
\beq
\det_{N}\big[ c_{\ell_a - m_b }[f]  \big] \cdot \Big( \det_{N}\big[ c_{a - b }[f]  \big]  \Big)^{-1}  \;. 
\enq

Although there should, quite probably, exist a direct transformation that would allow to connect Tracy-Widom's answer to ours
(and hence Bump and Diaconis one due to \cite{DehayeProofIdentityBumpDiaconisTracyWidomLacunatyToeplitz}), 
we have not succeeded in finding it. 
The formula can, of course, be simplified as soon as one provides some more informations on the lacunary integers 
$p_a, h_a, k_b, t_b$. Below, we treat a specific example of such a simplification, much in the spirit of 
the one outlined in Corollary \ref{Corolaire Toeplitz lacunaire à un indice}. 
In order to state the theorem, we first need to introduce a matrix $\mc{N}^{(\eps)}_{A ; B}\big( \mc{J} ; c \big)$
depending on the sets of integers
\beq
\mc{J} \; = \; \Big\{ \{p_a \} \; ; \; \{h_a \}  \; ; \; \{ k_a \} \; ; \; \{t_a\} \Big\} \;. 
\enq
The set $\mc{J}$ parametrizes its entries according to 
\bem
\Big( \mc{N}^{(\eps)}_{I ; I}\big( \mc{J} ; c \big) \Big)_{ab} \; = \; 
-\eps \Oint{  \Dp{}\mc{D}_{1} }{} \f{ 1 }{ z(1-\eps 0^+)-s } \Bigg\{ 
\de_{b \leq c} \de_{a \leq c} \bigg( \f{ \a_{-\eps} (z) }{  \a_{-\eps} (s) } \bigg)^{\eps} s^{\eps h_a-1 } z^{\eps h_b-1}
\, + \, \de_{b>c} \bigg( \f{ \a_{\eps} (z) }{  \a_{\eps} (s) } \bigg)^{\eps} s^{-\eps p_a } z^{-\eps h_b} \\
\, - \, \de_{b>c} \de_{a\leq c} \bigg( \f{ \a_{\eps} (z) }{  \a_{-\eps} (s) } \bigg)^{\eps} s^{\eps h_a-1 } z^{-\eps h_b}   \Bigg\} 
\cdot \f{ \dd s \cdot \dd z }{ (2i\pi)^2 } 
\; +   \; \eps \cdot \de_{ b \leq c} \Oint{  \Dp{}\mc{D}_{1} }{} \f{ s^{-\eps p_a } \cdot z^{\eps k_b-1} }{ z(1+\eps 0^+)-s } 
  \bigg( \f{ \a_{-\eps} (z) }{  \a_{\eps} (s) } \bigg)^{\eps}  \cdot \f{ \dd s \cdot \dd z }{ (2i\pi)^2 } 
\end{multline}
\bem
\Big( \mc{N}^{(\eps)}_{I ; II}\big( \mc{J} ; c \big) \Big)_{ab} \; = \; 
-\eps \Oint{  \Dp{}\mc{D}_{1} }{} \f{ 1 }{ z(1-\eps 0^+)-s } \Bigg\{ 
\de_{b \leq c} \de_{a \leq c} \bigg( \f{ \a_{\eps} (z) }{  \a_{-\eps} (s) } \bigg)^{\eps} s^{\eps h_a-1 } z^{- \eps t_b}
\, + \, \de_{b > c} \de_{a \leq c} \bigg( \f{ \a_{-\eps} (z) }{  \a_{-\eps} (s) } \bigg)^{\eps} s^{\eps h_a-1 } z^{\eps k_b-1} \\
\, - \, \de_{ b \leq c} \bigg( \f{ \a_{\eps} (z) }{  \a_{\eps} (s) } \bigg)^{\eps} s^{ -\eps p_a } z^{-\eps t_b}   \Bigg\} 
\cdot \f{ \dd s \cdot \dd z }{ (2i\pi)^2 } 
\; +   \; \eps \cdot \de_{ b > c} \Oint{ \Dp{}\mc{D}_{1} }{} \f{ s^{-\eps p_a } \cdot z^{\eps k_b-1 } }{ z(1+\eps 0^+)-s } 
  \bigg( \f{ \a_{-\eps} (z) }{  \a_{\eps} (s) } \bigg)^{\eps}  \cdot \f{ \dd s \cdot \dd z }{ (2i\pi)^2 } 
\end{multline}
and finally
\beq
\Big( \mc{N}^{(\eps)}_{II ; I}\big( \mc{J} ; c \big) \Big)_{ab} \; = \; 
- \eps \Oint{  \Dp{}\mc{D}_{1} }{} \f{ s^{\eps t_a-1 } }{ z(1-\eps 0^+)-s } \Bigg\{ 
\de_{b \leq c}  \bigg( \f{ \a_{-\eps} (z) }{  \a_{-\eps} (s) } \bigg)^{\eps}  z^{\eps k_b-1  }
\, - \, \de_{b > c}  \bigg( \f{ \a_{\eps} (z) }{  \a_{-\eps} (s) } \bigg)^{\eps}  z^{- \eps h_b } 
\Bigg\}\cdot \f{ \dd s \cdot \dd z }{ (2i\pi)^2 } 
\enq
\beq
\Big( \mc{N}^{(\eps)}_{II ; II}\big( \mc{J} ; c \big) \Big)_{ab} \; = \; 
- \eps \Oint{  \Dp{}\mc{D}_{1} }{} \f{ s^{\eps(t_a-1)} }{ z(1-\eps 0^+)-s } \Bigg\{ 
\de_{b \leq c}  \cdot \bigg( \f{ \a_{\eps} (z) }{  \a_{-\eps} (s) } \bigg)^{\eps} \cdot  z^{- \eps t_b }
\, + \, \de_{b > c} \cdot  \bigg( \f{ \a_{-\eps} (z) }{  \a_{-\eps} (s) } \bigg)^{\eps} \cdot  z^{ \eps k_b-1 } 
\Bigg\}\cdot \f{ \dd s \cdot \dd z }{ (2i\pi)^2 } 
\enq

The $(1-\eps O^{+})$ prescription means that the integral should be understood as the limit when $z$ approaches a point 
on $\Dp{}\mc{D}_{1}$ from the inside ($\eps=+1$) or outside $(\eps=-1)$ of the unit disk. 

\begin{theorem}
Assume that the lacunary integers $p_a,h_a$ are given as in \eqref{particules trous pour entiers pa} and \eqref{particuler trous pour entiers ha}
and, likewise, that the lacunary integers $k_b, t_b$ are given as :
\beqa
k_a   \; = \; 1-k_a^{-} \quad \e{for} \quad a=1,\dots, r_{-} \qquad & \e{and} & 
\qquad k_{a+r_-} \; = \;  k_a^{+} + N  \quad \e{for} \quad a=1,\dots, r_{+}  
\label{particules trous pour entiers pa} \\
t_a   \; = \; t_a^{-}  \quad \e{for} \quad a=1,\dots, r_{-} \qquad &\e{and} &
				\qquad t_{a+r_-} \; = \; N+1- t_a^{+}  \quad \e{for} \quad a=1,\dots, r_{+}   \;.  
\label{particuler trous pour entiers ha}
\eeqa
Further, let the sets $\{h_a^{+}\}_1^{n_+}$ and $\{k_a^+\}_{1}^{r^+}$ (resp. $\{h_a^{-}\}_1^{n_-} $ and  $\{k_a^-\}_{1}^{r^-}$)
be well ordered with overlap $c_+$ (resp. $c_-$). Then, provided that the matrices 
$\mc{N}^{(\pm)}\big( \mc{J}^{(\pm)} ; c_{\pm} \big)$ have maximal rank,
one has the asymptotic expansion
\beq
\det_{N}\big[ c_{\ell_a - m_b }[f]  \big] \cdot \Big( \det_{N}\big[ c_{a - b }[f] \big]  \Big)^{-1} \; = \; 
\det_{n^+ + r^+} \big[  \mc{N}^{(+)}\big( \mc{J}^{(+)} ; c_{+} \big) \big] 
\cdot \det_{n^- + r^-} \big[  \mc{N}^{(-)}\big( \mc{J}^{(-)} ; c_{-} \big) \big]  \cdot 
\Big( 1 + \e{O}\big( N^{-\infty} \big) \Big) \;. 
\enq
The sets $\mc{J}^{(\pm)}$ appearing above are defined as
\beq
\mc{J}^{(+)} \; = \; \Big\{ \{p_a^+ \}_1^{n_+} \; ; \; \{h_a^+ \}_1^{n_+}  \; ; \; \{ k_a^+ \}_1^{r_+} \; ; \; \{t_a^+\}_1^{r_+} \Big\} 
\qquad  \e{and} \qquad
\mc{J}^{(-)} \; = \; \Big\{ \{1-p_a^- \}_1^{n_-} \; ; \; \{1-h_a^- \}_1^{n_-}  
		\; ; \; \{1- k_a^- \}_1^{r_-} \; ; \; \{1-t_a^-\}_1^{r_-}   \Big\}  \;. 
\nonumber
\enq
\end{theorem}

\Proof 

The representation obtained in proposition \ref{Proposition factorisation Toeplitz lacunaire a ligne et colonnes} 
is invariant under a permutation of the integers $h_a, t_a$ \textit{along} with a simultaneous permutation 
of the associated integers $p_a, k_a$. Hence, reorganising the entries of the matrix $\mc{N}$ in each block 
so that the natural order imposed by the $\pm$ splitting of the lacunary integers is respected and a 
repeated application of the manipulations outlined in the proof of Theorem \ref{Theorem Asymptotiques Toeplitz lacunaire a lignes}
leads to $\det_{n+r}\big[ \mc{N} \big]  \; = \;   \det_{n+r}\big[ \wh{\mc{N} } \big]  $ where 
\beq
 \wh{\mc{N} } \; = \;  \left(   \ba{c c}  
\ba{cc}  \mc{N}_{I ; I}^{(-)}\Big(\mc{J}^{(-)} ; c_- \Big)   &   0 \\
									0    & \mc{N}_{I ; I}^{(+)}\Big(\mc{J}^{(+)} ; c_+ \Big)    \ea 
& \ba{cc}  \mc{N}_{I ; II}^{(-)}\Big(\mc{J}^{(-)} ; c_- \Big)   &  0 \\
										0   & \mc{N}_{I ; II}^{(+)}\Big(\mc{J}^{(+)} ; c_+ \Big)   \ea  \\ 
\ba{cc}  \mc{N}_{II ; I}^{(-)}\Big(\mc{J}^{(-)} ; c_- \Big)   &   0 \\
									0    & \mc{N}_{II ; I}^{(+)}\Big(\mc{J}^{(+)} ; c_+ \Big)    \ea 
& \ba{cc}  \mc{N}_{II ; II}^{(-)}\Big(\mc{J}^{(-)} ; c_- \Big)   &  0 \\
								0   & \mc{N}_{II ; II}^{(+)}\Big(\mc{J}^{(+)} ; c_+ \Big)   \ea 	\ea	\right)
\;  + \;   \e{O}\big(N^{-\infty}\big) \;. 
\enq
Finally, $ \e{O}\big(N^{-\infty}\big) $ appearing above refers to an $(n+r)\times (n+r)$ matrix whose all 
entries are a $ \e{O}\big(N^{-\infty}\big) $. It is then enough to exchange the appropriate lines and columns
in $\det_{n+r}\big[ \wh{\mc{N} } \big]  $ and invoke the maximality of the rank of the matrices $\mc{N}^{(\pm)}$. \qed


\section*{Conclusion}

In this paper we have proposed a Riemann--Hilbert problem based approach to the analysis of the large-size asymptotic
behaviour of lacunary Topelitz determinants having a \textit{finite} number of modified lines an rows. 
Our approach allows one to obtain an alternative to the ones obtained in 
\cite{BumpDiaconisLacunaryToeplitzThrougSumsSymFctsAndYoungTableaux,TracyWidomAsymptoticExpansionLacunaryToeplitz} 
representation for its large-$N$ asymptotics. 
Our answer involves a determinant that solely depends on the number of modified rows and lines and not 
on the index of the largest modified line or column. In particular, this allows one to investigate the asymptotics in the 
case when the locii of some of the modified lines and columns go to infinity. We have treated certain instances 
of such a situation in the present paper. It is clear from the very setting of our analysis that our method
allows one to treat also generalisations of lacunary Toeplitz determinants such as those considered in 
\cite{LionsToeplitzLacunaires}. To do so, one should simply replace the functions 
$\wh{\tau}_{t_a}$ and $\wh{\kappa}_{h_a}$ arising in our analysis by more general ones.

\section*{Acknowledgements}

The author is supported by CNRS and acknowledges support from the 
the Burgundy region PARI 2013 FABER grant "Structures et asymptotiques d'int\'{e}grales multiples".
The author would like to thank J.M.-Maillet and M. Piatek 
for stimulating discussions. The author is indebted to the Joint Institute for 
Nuclear Research, Dubna, for its warm hospitatlity during the time when part of this
work has been done.






\appendix

\section{Asymptotic inversion of $I+V_0$-The Riemann--Hilbert approach}
\label{Appendix RHP pour Toeplitz regulier}

\subsection{The Riemann--Hilbert problem associated with $I+V_0$}

Consider the Riemann--Hilbert problem for a piecewise analytic $2\times 2$ matrix $\chi$ having a jump on the unit circle 
$\Dp{}\mc{D}_1$:
\begin{itemize}
\item $\chi$ is analytic on   $\Cx\setminus \Dp{}\mc{D}_1$ \; ;
\item $\chi(z) \; = \;  I_2  \; + \; \f{1}{z} \cdot \e{O}\bigg( \ba{cc} 1 & 1 \\ 1 & 1\ea  \bigg) \; $ when $z \tend \infty $;
\item $\chi$ admits continuous $\pm$-boundary values on $\Dp{}\mc{D}_1 $;
\item $\chi_{+}(z) \cdot  \left(\ba{cc}
                   2-f(z) & \big( f(z)-1 \big)\cdot  z^N \\
                   \big( 1 - f(z) \big)\cdot  z^{-N} & f(z)
                   \ea \right) \; = \; \chi_-(z)  \;\; ; \quad z \in \Dp{}\mc{D}_1 \, .$
\end{itemize}
In the formulation of the Riemann--Hilbert problem, we have adopted the following notations. 
Given an oriented Jordan curve $\Ga \subset \Cx$ and a function $f$
on $\Cx\setminus \Ga$, $f_{\pm}$ refer to the $\pm$-boundary values of the function $f$ on $\Ga$ where $+$ (resp. $-$)
refers to approaching a point on $\Ga$ non-tangentially  from the left (resp. right) side of the curve. 
Finally, a matrix domination of the sort $A = \e{O}(B)$ is to be understood entry-wise \textit{viz}
$A_{jk}=\e{O}(B_{jk})$. 

We also remind that the unit circle $\Dp{}\mc{D}_1$ is oriented canonically (\textit{ie} the $+$ side of the
contour corresponds to the interior of the circle). It is a standard fact that 
the above Riemann--Hilbert problem admits a unique solution. 

\subsection{Transformation to a perturbatively solvable Riemann--Hilbert problem for $\Ups$}

We now define a new matrix $\Upsilon$ according to Fig.~\ref{Contour du RHP pour Y}, \textit{ie}
\begin{itemize}
\item $\Upsilon=\chi \, \a^{\sg_3}\, , $ for $z$ being in the exterior of $\Ga_{\e{ext}}$ and the interior of $\Ga_{\e{int}}\; ;$
\item $\Upsilon=\chi \, \a^{\sg_3} M_{\e{ext} }^{-1}\, , $ for $z$  between $\Ga_{\e{ext}}$ and $ \Dp{}\mc{D}_{1}\; ;$
\item $\Upsilon=\chi \, \a^{\sg_3} M_{ \e{int}}\, , $ for $z$ between $\Ga_{\e{int}}$ and $ \Dp{}\mc{D}_{1} $.
\end{itemize}
Here $\a$ is as defined in \eqref{definition facteur alpha solution RHP scalaire}. It is readily seen that it
solves the scalar RHP
\beq
\a\;  \e{analytic}\; \e{on} \; \Cx\setminus \Dp{}\mc{D}_1 \quad \a_- = f  \a_+  \; ,\qquad \e{on} \; \; \Dp{}\mc{D}_1
\quad \a (z) \tend 1 \; \e{when}\; z \tend \infty \;.
\enq
  The matrices
 $M_{\e{int}/\e{ext}}$ appearing in the definition of $\Ups$ read
\beq
 M_{ \e{int} }(z) \; = \;  \left(\ba{cc}
                         			1& \big( 1-f^{-1}(z) \big)\a^{-2}(z) \cdot  z^N \\
					                         0& 1 \ea\!\! \right), \quad
 M_{\e{ext}}(z) \; = \; \left(\ba{cc}
                        1   &    0\\
                         \big( f^{-1}(z) - 1) \a^{2}(z)\cdot  z^{-N} &  1 \ea \!\!\right)\; .
\enq

The curves $\Ga_{\e{ext}}$ and $\Ga_{\e{int}}$ are chosen in such a way that they are located inside of the open neighbourhood 
of $\Dp{}\mc{D}_1$ on which $f$ is holomorphic.
\begin{figure}[h]
\centering
\begin{pspicture}(9,8)

\pscircle(4,4){4}
\psline[linewidth=3pt]{->}(4,8)(3.95,8)

\pscircle[linestyle=dashed](4,4){2.5}
\psline[linewidth=3pt]{->}(6.5,4)(6.5,4.1)

\pscircle(4,4){1}
\psline[linewidth=3pt]{->}(4,5)(3.95,5)

\rput(4,4){$\Upsilon=\chi \cdot \a^{\sg_3}$}

\rput(4,7){$\Upsilon=\chi \cdot \a^{\sg_3} \cdot  M_{ \e{ext} }^{-1}$}

\rput(4,2.5){$\Upsilon=\chi \cdot \a^{\sg_3} \cdot M_{ \e{int} } $}

\rput(7.5,7.5){$\Upsilon=\chi \cdot \a^{\sg_3}$}


\rput(8.4,4.3){$\Ga_{\e{ext}}$}

\rput(5.3,4){$\Ga_{ \e{int} } $}

\rput(1.3,4){$\Dp{}\mc{D}_1$}

\end{pspicture}
\caption{ Contour for the RHP $\Upsilon$ and the associated contour $\Ga_{\Ups}=\Ga_{\e{int}}\cup\Ga_{\e{ext}}$. 
\label{Contour du RHP pour Y}}
\end{figure}
One readily sees that $\Upsilon$ satisfies the RHP
\begin{itemize}
\item $\Upsilon$ is analytic in $\Cx\setminus \Ga_{\Ups}\; $;
\item $\Ups(z) \; = \;  I_2  \; + \; \f{1}{z} \cdot \e{O}\bigg( \ba{cc} 1 & 1 \\ 1 & 1\ea  \bigg) \; $ when $z \tend \infty $;
\item $\Ups$ admits continuous $\pm$-boundary values on $\Ga_{\Ups}$;
\item $\Ups_{+}(z) \cdot  G_{\Ups}(z) \; = \; \Ups_-(z)  \;,  \quad z \in \Ga_{\Ups} $ \; \; 
where \; \; 
$G_{\Ups}(z) \; = \; M_{ \e{ext} }(z) \cdot \bs{1}_{ \Ga_{\e{ext}} }(z) \; + \;  M_{ \e{int} }(z)\cdot \bs{1}_{ \Ga_{\e{int}} }(z) $.
\end{itemize}
and $\bs{1}_A$ stands for the indicator function of the set $A$. 
Since
\beq
\norm{ G_{\Ups}\, - \, I_2 }_{ L^{\infty}\big(\Ga_{\Ups}\big) } \; + \; 
\norm{ G_{\Ups}\, - \, I_2 }_{ L^{1}\big(\Ga_{\Ups}\big) } \; + \; 
\norm{ G_{\Ups}\, - \, I_2 }_{ L^{2}\big(\Ga_{\Ups}\big) } \; \leq \; C_1 \ex{- \kappa N}
\enq
for some constants $C_1>0$ and $\kappa>0$, it follows from the equivalence of the Riemann--Hilbert problem for $\Ups$ with the singular integral equation 
satisfied by $\Ups_{+}$ that, for any compact $K\supset \Ga_{\Ups}$, 
\beq
\norm{ \Ups \, - \, I_2 }_{ L^{\infty}\big(\Cx\setminus K \big) }  \; \leq \; C_1^{\prime} \ex{- \kappa N}
\enq
 for some new constant $C_1^{\prime} >0$.

\subsection{The resolvent operator and factorisation of the determinant}

It is well known \cite{ItsIzerginKorepinSlavnovDifferentialeqnsforCorrelationfunctions} that the resolvent 
kernel $R_0$ associated with the integrable integral operator $I+V_0$ takes the form 
\beq
R_0(z,s) \; = \; \f{ \big(\bs{F}_L (z) , \bs{F}_{R}(z) \big)   }{ z-s }
\enq
where given vector $\bs{x},\bs{y}$, $(\bs{x},\bs{y})=x_1y_1+x_2y_2$ and 
\beq
\bs{F}_L^{\bs{T}}(z) \; = \;  \bs{E}_L^{\bs{T}}(z)\cdot \chi^{-1}(z) \qquad 
\bs{F}_R(z) \; = \;  \chi(z )\cdot \bs{E}_R(z)
\enq
where the two-dimensional vectors $\bs{E}_R(z), \bs{E}_L(z)$ take the form
\beq
\bs{E}_L^{\bs{T}}(z) \; = \; \big( f(z)-1 \big)\cdot \Big( - z^{-\f{N}{2}} \, , \, z^{\f{N}{2}} \Big)
\qquad \e{and} \qquad 
\bs{E}_R^{\bs{T}}(z) \; = \; \f{ 1 }{ 2i\pi }\cdot \Big(  z^{\f{N}{2}} \, , \, z^{-\f{N}{2}} \Big) \;. 
\enq
It follows from the factorisation of $\chi$ that $\bs{F}_R$ can be recast as 
\beq
\bs{F}_R(z) \; = \; \bs{F}_R^{(0)}(z) \; + \; \bs{F}^{(\infty)}_R(z)  \qquad \e{with} \quad 
\left\{ \ba{ccc} \bs{F}_R^{(0)}(z) &= &  M_{\e{int};+}^{-1}(z) \cdot  \a_+^{-\sg_3}(z) \cdot \bs{E}_R(z) \vspace{2mm}\\
	\bs{F}_R^{(\infty)}(z) &= &  \big( \Ups-I_2 \big) \cdot M_{\e{int};+}^{-1}(z) \cdot \a_+^{-\sg_3}(z) \cdot  \bs{E}_R(z) \ea \right. \:. 
\enq

The uniform bounds on $\Ups - I_2$ ensure that 
\beq
\big| \big| \bs{F}_R^{(\infty)}  \big|\big|_{L^{\infty}(\msc{C}) } \; = \; C\cdot\ex{-\kappa N}
\enq
 for some $C>0$. 
Thus, for $z \in \msc{C}$, 
\beq
\bs{F}_R^{(0)}(z) \; = \; \f{1}{2i\pi} \cdot \left( \ba{c} z^{\f{N}{2}} \a_{-}^{-1}(z)  \\ z^{-\f{N}{2}} \a_{+}(z)  \ea \right)
\quad \e{and} \quad
\bs{E}_L^{(0)}(z) \; = \; \big( b(z) -1 \big) \cdot \left( \ba{c} - z^{-\f{N}{2}} \a_{+}(z)  \\ z^{ \f{N}{2} } \a_{-}^{-1}(z)  \ea \right)
\enq

As a consequence, the resolvent kernel $R_0$ decomposes exactly as given in 
\eqref{ecriture decomposition resolvant en partie finie et perturbation exp petite en N}.

\end{document}